\documentclass[12pt, onesided, reqno]{amsart}

\usepackage{amssymb}
\usepackage{amsmath}
\usepackage{color}

\evensidemargin\oddsidemargin

\begin{document}

\setcounter{page}{1}

\newtheorem{PROP}{Proposition}
\newtheorem{REMS}{Remark}
\newtheorem{LEM}{Lemma}
\newtheorem{THEA}{Theorem A\!\!}
\renewcommand{\theTHEA}{}
\newtheorem{THEB}{Theorem B\!\!}
\renewcommand{\theTHEB}{}

\newtheorem{theorem}{Theorem}
\newtheorem{proposition}[theorem]{Proposition}
\newtheorem{corollary}[theorem]{Corollary}
\newtheorem{lemma}[theorem]{Lemma}
\newtheorem{assumption}[theorem]{Assumption}

\newtheorem{definition}[theorem]{Definition}
\newtheorem{hypothesis}[theorem]{Hypothesis}

\theoremstyle{definition}
\newtheorem{example}[theorem]{Example}
\newtheorem{remark}[theorem]{Remark}
\newtheorem{question}[theorem]{Question}

\newcommand{\eqnsection}{
\renewcommand{\theequation}{\thesection.\arabic{equation}}
    \makeatletter
    \csname  @addtoreset\endcsname{equation}{section}
    \makeatother}
\eqnsection

\def\a{\alpha}
\def\b{\beta}
\def\B{{\bf B}} 
\def\cia{c_{\a, \infty}}
\def\coa{c_{\a, 0}}
\def\cua{c_{\a, u}}
\def\cL{{\mathcal{L}}} 
\def\Ea{E_\a}
\def\eps{\varepsilon}
\def\g{{\gamma}} 
\def\Ga{{\Gamma}} 
\def\i{{\rm i}}
\def\K{{\bf K}}
\def\Ka{{\bf K}_\a}
\def\L{{\bf L}}
\def\lbd{\lambda}
\def\lcr{\left[}
\def\lpa{\left(}
\def\lva{\left|}
\def\rpa{\right)}
\def\rcr{\right]}
\def\rva{\right|}
\def\T{{\bf T}}
\def\M{{\mathcal M}}
\def\X{{\bf X}}
\def\U{{\bf U}}
\def\Un{{\bf 1}}
\def\ZZ{{\bf Z}}
\def\CC{{\bf C}}
\def\GG{{\bf \Ga}}
\def\BB{{\bf B}}
\def\car{c_{\a,\rho}}
\def\sar{s_{\a,\rho}}
\def\pbxy{\Pb_{(x,y)}}
\def\foxy{f^0_{x,y}}

\def\E{\mathbb{E}}
\def\Z{\mathbb{Z}}
\def\N{\mathbb{N}}
\def\Q{\mathbb{Q}}
\def\R{\mathbb{R}}
\def\Pb{\mathbb{P}}
\def\C{\mathbb{C}}
\def\F{\mathcal{F}}
\def\S{\mathcal{S}}
\def\W{\mathcal{W}}
\def\L{\mathcal{L}}
\def\G{\mathcal{G}}

\newcommand{\equi}{\mathop{\sim}\limits}
\def\={{\;\mathop{=}\limits^{\text{(law)}}\;}}
\def\d{{\;\mathop{=}\limits^{\text{(1.d)}}\;}}
\def\st{{\;\mathop{\geq}\limits^{\text{(st)}}\;}}

\def\claw{\stackrel{d}{\longrightarrow}}
\def\elaw{\stackrel{d}{=}}
\def\qed{\hfill$\square$}
                  
\title[Persistence of integrated stable processes]
      {Persistence of integrated stable processes}

\author[Christophe Profeta]{Christophe Profeta}

\address{Laboratoire d'analyse et probabilit\'es, Universit\'e d'Evry-Val d'Essonne, B\^atiment IBGBI, 23 boulevard de France, F-91037 Evry Cedex. {\em Email} : {\tt christophe.profeta@univ-evry.fr}}

\author[Thomas Simon]{Thomas Simon}

\address{Laboratoire Paul Painlev\'e, Universit\'e Lille 1, F-59655 Villeneuve d'Ascq Cedex and Laboratoire de physique th\'eorique et mod\`eles statistiques, Universit\'e  Paris-Sud, F-91405 Orsay Cedex. {\em Email} : {\tt simon@math.univ-lille1.fr}}



\begin{abstract} We compute the persistence exponent of the integral of a stable L\'evy process in terms of its self-similarity and positivity parameters. This solves a problem raised by Z. Shi (2003). Along the way, we investigate the law of the stable process $L$ evaluated at the first time its integral $X$ hits zero, when the bivariate process $(X,L)$ starts from a coordinate axis. This extends classical formul\ae\,\! by McKean (1963) and Gor'kov (1975) for integrated Brownian motion.

\end{abstract}

\maketitle
 
\section{Introduction and statement of the results}

Let $X =\{ X_t, \, t\ge 0\}$ be a real process starting at zero and $T_x = \inf\{ t > 0, \; X_t > x\}$ be its first-passage time above a positive level $x$. Studying the law of $T_x$ is a classical problem in probability theory. In general, it is difficult to obtain an explicit expression of this law. However, it has been observed that in many interesting cases the survival function has a polynomial decay:
\begin{equation}
\label{Perse}
\Pb[T_x > t]\; =\; t^{-\theta + o(1)}, \qquad t \to +\infty,
\end{equation}
where $\theta$ is a positive constant which is called the persistence exponent, and usually does not depend on $x$. The computation of persistence exponents has many connections to various problems in probability and mathematical physics, and we refer to the recent surveys \cite{AS, BMS} for more information on this topic. In this paper we consider this problem for the process
$$X_t \; =\; \int_0^t L_s\, ds,$$
where $L = \{L_t, \, t\ge 0\}$ is a strictly $\a-$stable L\'evy process starting from zero, with law $\Pb.$ Our process $L$ is normalized to have characteristic exponent
\begin{equation}
\label{Norm}
\Psi(\lbd)\; =\;\log(\E[e^{\i \lbd L_1}])\; =\; -(\i \lbd)^\a e^{-\i\pi\a\rho\, {\rm sgn}(\lbd)}, \qquad \lbd\in\R,
\end{equation}
where $\a\in (0,2]$ is the self-similarity parameter and $\rho = \Pb[L_1 \ge 0]$ is the positivity parameter. We refer to \cite{ST, Z} for classic accounts on stable laws and processes. The strict stability implies the $(1/\a)-$self-similarity of $L$ and the $(1+1/\a)-$self-similarity of $X$, in other words that
$$\{L_{kt}, \; t\ge 0\}\; \elaw\; \{k^{1/\a}L_{t}, \; t\ge 0\}\qquad\mbox{and}\qquad \{X_{kt}, \; t\ge 0\}\; \elaw\; \{k^{1+1/\a}X_{t}, \; t\ge 0\}$$
for all $k > 0.$ When $\a = 2,$ one has $\rho = 1/2$ and $\Psi(\lbd) = -\lbd^2,$ so that $L = \sqrt{2} B$ is a rescaled Brownian motion. When $\a = 1,$ one has $\rho\in (0,1)$ and $L$ is a Cauchy process with a linear drift. When $\a\in (0,1)\cup(1,2)$ the characteristic exponent takes the more familiar form
$$\Psi(\lbd) \; =\; -\kappa_{\a,\rho}\vert\lbd\vert^\a (1 - \i\b\tan(\pi\a/2)\,{\rm sgn}(\lbd)),$$
where $\b\in[-1,1]$ is an asymmetry parameter, whose connection with the positivity parameter is given by Zolotarev's formula:
$$\rho \; =\; \frac{1}{2} \,+ \,\frac{1}{\pi\a} \arctan(\b\tan(\pi\a/2)),$$
and $\kappa_{\a,\rho} = \cos(\pi\a(\rho -1/2)) > 0$ is a scaling constant. The latter could have taken any positive value, changing the normalization (\ref{Norm}) accordingly, without incidence on our purposes below. One has $\rho \in [0,1]$ if $\a < 1$ and $\rho\in[1-1/\a, 1/\a]$ if $\a > 1.$ When $\a > 1$ and $\rho = 1/\a$ the process $L$ has no positive jumps, whereas it has no negative jumps when $\a > 1$ and $\rho = 1-1/\a$. When $\a < 1$ and $\rho = 0$ or $\rho = 1,$ the process $\vert L\vert$ is a stable subordinator and has increasing sample paths, a situation which will be implicitly excluded throughout this paper. In this case, the process $X$ is indeed also monotonous and the survival function in (\ref{Perse}) either is one or decays towards zero at an exponential speed - see \cite{AS} p.4 for details. \\
 
When $\a =2,$ the bivariate process $(X,L)$ is Gaussian with explicit covariance function and transition density, providing also the basic example of a degenerate diffusion process - see \cite{Lab} for details and references. When $\a < 2,$ the process $(X,L)$ is Non-Gaussian $\a-$stable in the broad sense of \cite{ST}. 
The process $(X,L)$ is a strong Markov process, which is sometimes called the 
Kolmogorov process in the literature. In the following we will set $\pbxy$ for the law of $(X,L)$ starting at $(x,y)\in\R^2.$ Our main concern in this paper is the hitting time of zero for $X:$ 
$$T_0\; =\; \inf\{ t> 0, \; X_t = 0\}.$$
Since $\vert L\vert$ is not a subordinator, a simple argument using self-similarity and the zero-one law for Markov processes - see below Lemma \ref{Klar} for details - shows that $\Pb_{(0,0)}[T_0 = 0] = 1,$ in other words that the origin is regular for the vertical axis. If $x < 0$ or $x=0$ and $y <0,$ the continuity of the sample paths of $X$ show that a.s. $T_0 = \inf\{ t> 0, \; X_t > 0\},$ and it will be checked in Lemma \ref{Klar} below that $T_0$ is also a.s. finite. If $x > 0$ or $x=0$ and $y >0,$ the law of $T_0$ is obviously deduced from that of the latter situation in considering the dual L\'evy process $-L.$ 

When $(x,y)\neq (0,0),$ the difficulty to obtain concrete informations on the law of $T_0$ under $\pbxy$ comes from the fact that $X$ itself is  not a Markov process. In the Brownian case for example, the density function of $T_0$ is expressed through quite intricate integral formul\ae\, - see \cite{AS} pp.15-16 and the references therein. On the other hand, some universal estimates can be obtained for the behaviour of the distribution function $\pbxy[T_0\le t]$ as $t\to 0,$ using self-similarity and Gaussian or stable upper tails for the supremum process - see e.g. Section 10.4 in \cite{ST}. But it is well-known that the study of $\pbxy[T_0 > t]$ as $t\to +\infty$ is a harder problem, where a more exotic behaviour is expected.\\

Throughout the paper, for any real functions $f$ and $g$ we will use the standard notation $f(t)\asymp g(t)$ as $t\to +\infty$ to express the fact that there exist $\kappa_1, \kappa_2$  positive finite constants such that $\kappa_1 f(t)\le g(t) \le \kappa_2 f(t)$ as $t\to +\infty.$ Our main result is the following. 

\begin{THEA} Assume $x < 0$ or $x=0,y < 0.$ One has
$$\Pb_{(x,y)} [T_0 > t] \; \asymp\; t^{-\theta},\qquad t\to +\infty,$$
with $\theta = \rho/(1+ \a(1-\rho)).$ 
\end{THEA}

In the Brownian case $\a = 2,$ one has $\theta = 1/4=\rho/2$ and this estimate has been known since the works of M.~Goldman - see Proposition 2 in \cite{G}, with a more precise formulation on the density function of $T_0$, following the seminal article of McKean \cite{McK}. The universality of the persistence exponent 1/4 for integrals of real L\'evy processes having exponential moments on both sides has been shown in \cite{AD}, with the help of strong approximation arguments. Recently, it was proved in \cite{DDG} that all integrated real random walks with finite variance have also 1/4 as persistence  exponent, extending \cite{Si} for the particular case of the integrated simple random walk. Let us also mention that the survival function of the $n^{\rm th}$ hitting time of zero for the integrated Brownian motion exhibits the same power decay up to a logarithmic term in $ct^{-1/4}(\ln(t))^{n-1}$ with an explicit constant $c$, as shown by the first author in \cite{CP1}.

In the case $1 < \a < 2$ and with no negative jumps, that is $\rho = 1 -1/\a,$ one obtains $\theta = (\a - 1)/2\a = \rho/2,$ an estimate which had been proved by the second author in \cite{TS1} with different techniques and a less precise formulation than Theorem A for the lower bound, involving a logarithmic correction term. It is worth mentioning that the same persistence exponent $(\a - 1)/2\a$ appears for the integrals of  random walks attracted towards this spectrally positive L\'evy process - see Remark 1.2 in \cite{DDG}.\\

It has been conjectured in \cite{AS} - see Conjecture 4 therein - that the persistence exponent should be $\rho/2$ in general. This expected value should be compared with a classical result of Bingham stating that the persistence exponent is $\rho$ for the stable process $L$ - see (2.16) in \cite{AS} and the references of Section 2.2 therein. The admissible set of $(\a,\rho)$ and Theorem A entail that $\theta > \rho/2$ as soon as $L$ has negative jumps, hence providing a negative answer to this  conjecture. The fact that $\theta$ is an increasing function of the positivity parameter $\rho$ matches the intuition, however it is harder to explain heuristically why it is also a decreasing function of $\a.$ 

Specifying $x = -1$ and $y = 0$ in Theorem A entails by self-similarity the following lower tail probability estimate
$$\Pb [ X_1^*\le \eps]\; \asymp\; \eps^{\frac{\theta\a}{\a +1}}, \qquad \eps \to 0,$$
with the notation $X_1^* = \sup\{X_t, \, t\le 1\}.$ Some heuristics on the subordination of $X$ by the inverse local time of $L$ when $\a > 1$ had led to the conjecture, formulated in Part 2 of \cite{ZS}, that in the symmetric case $\rho =1/2$ one should have $\Pb [ X_1^*\le \eps] = \eps^{(\a-1)^+/2(\a +1) +o(1)}$ as $\eps\to 0.$ The invalidity of this conjecture as soon as $\a$ is close enough to 1 had been observed in \cite{TS2}. Theorem A shows that Shi's exponent is the right one only for integrated Brownian motion: in the symmetric case one has $\theta\a/(\a +1) = \a/(\a+1)(\a+2) \ge (\a-1)^+/2(\a+1),$ with an equality only if $\a =2.$ Let us mention in passing that lower tail probabilities offer some challenging problems for Gaussian processes - see \cite{Ma, QS}.\\
  
Our method to prove Theorem A hinges upon the random variable $L_{T_0},$ the so-called hitting place of $(X,L)$ on the vertical axis, which has been extensively studied in the Brownian case - see \cite{McK,G, La, Lab}. Notice that this random variable is positive under $\pbxy$ if $x<0$ or $x=0$ and $y <0$. The reason why it is connected to the persistence exponent comes from the following heuristical equivalence for fractional moments
$$\E_{(x,y)}[T_0^s]\, <\, +\infty\quad \Leftrightarrow\quad \E_{(x,y)}[L_{T_0}^{\a s}]\, <\, +\infty$$
for all $s >0,$ which had been conjectured in \cite{TS1} p.176, and turns out to be true as a consequence of Theorem A and Lemma \ref{lem:asympLT0} below. The precise relationship between the upper tails of $T_0$ and that of $L_{T_0}$ follows from a series of probabilistic estimates which are the matter of Section 4.

In this paper we also provide a rather complete description of the law of the random variable $L_{T_0}$ when $(X,L)$ starts from a coordinate axis. To express our second main result, we need some further notation. For every $\mu\in (0,1),$ introduce the $\mu-$Cauchy random variable $\CC_\mu,$ with density
$$\frac{\sin(\pi\mu)}{\pi\mu (x^2 + 2 \cos (\pi\mu) x + 1)}1_{\{x\geq0\}}.$$
Our above denomination comes from the case $\mu = 1/2,$ where $\CC_{1/2}$ is the half-Cauchy distribution. If $X$ is a positive random variable and $\nu\in\R$ is such that $\E[X^\nu] < \infty,$  the positive random variable $X^{(\nu)}$ defined by
$$\E[f(X^{(\nu)})]\; =\; \frac{\E[X^\nu f(X)]}{\E[X^\nu]}$$
for all $f : \R^+ \to \R$ bounded continuous, is known as the size bias of order $\nu$ of $X.$ Observe that when $X$ is absolutely continuous, the density of $X^{(\nu)}$ is obtained in multiplying that of $X$ by $x^\nu$ and renormalizing. Introduce finally the parameters 
$$\g \, =\, \frac{\rho\a}{1+\a}\qquad\mbox{and}\qquad \chi\, =\, \frac{\rho\a}{1+\a(1-\rho)}=\alpha\theta.$$
Notice that from the admissible set for $(\a, \rho),$ we have $\g \in (0,1/2)$ and $\chi\in (0,1).$ 

\begin{THEB} {\em (i)} For every $y < 0,$ under $\Pb_{(0,y)}$ one has
$$L_{T_0}\, \elaw\, \vert y\vert (\CC_\chi^{1-\g})^{(1)}.$$

{\em (ii)} For every $x < 0,$ under $\Pb_{(x,0)}$ the positive random variable $L_{T_0}$ has Mellin transform
$$\E_{(x,0)}[L_{T_0}^{s-1}]\; =\;  \frac{(1+\alpha)^{\frac{1-s}{1+\alpha}}\Ga(\frac{\a +2}{\a +1})\Ga (\frac{1-s}{\a +1})\sin (\pi\g)}{\Ga(\frac{s}{\a +1})\Ga (1-s)\sin (\pi s(1-\g))}\,\vert x\vert^{\frac{s-1}{\a+1}}, \qquad \vert s\vert < 1/(1-\g).$$
\end{THEB}

The proof of this result is given in Section 3, following some preliminary computations involving oscillating integrals and the Fourier transform of $X_t,$ performed in Section 2. Observe that the density in (i) above is explicit and reads for example
$$\frac{3}{2\pi} \frac{\vert y\vert^{1/2} z^{3/2}}{\vert y\vert^3 + z^3}1_{\{z\geq0\}}$$
in the Brownian case, a formula originally proved by McKean in \cite{McK} - see also formul\ae\, (1) and (2) in \cite{La}. As is well-known, the Cauchy random variable appears in exit or winding problems for two-dimensional Brownian motion. The fact that it is also connected with similar problems for general integrated stable processes is perhaps more surprising. 

An interesting consequence of (ii) is that the Mellin transform can be inverted in the Cauchy case $\a = 1$ and exhibits the same type of law as in (i): one obtains
$$L_{T_0}\, \elaw\, \sqrt{2\vert x\vert}\, (\CC_\delta^{1-\g})^{(1)}$$
with the notation $\delta = (1+\chi)/2.$ The Mellin transform of (ii)  can also be simply inverted in the Brownian case in terms of Beta and Gamma random variables, shedding some new light on a formula by Gor'kov \cite{Gk} which was of the analytical type, and in the case $\a < 1$ in terms of  positive stable random variables. The Mellin inversion is however more complicated when $\a\in (1,2),$ and involves no classical random variables in general - see Section 3.3 below for details.

\section{Preliminary computations}

The following lemma, which we could not locate in the literature, will be useful in the sequel.
\begin{lemma} 
\label{switch}
Let $\nu\in (0,1)$ and $X$ be a real random variable such that $\E[\vert X\vert^{-\nu}] < \infty.$ One has
$$\int_0^{\infty}  \lbd^{\nu-1} \,\E[\cos(\lbd X)] \, d\lbd\; =\; \Ga(\nu) \cos (\pi\nu/2)\,\E[ \vert X\vert^{-\nu}]$$
and
$$\int_0^{\infty}  \lbd^{\nu-1} \E[\sin(\lbd X)]\, d\lbd\; =\; \Ga(\nu) \sin (\pi\nu/2)\,\E[ \vert X\vert^{-\nu}{\rm sgn}(X)].$$
\end{lemma}

\proof
The generalized Fresnel integral which is computed e.g. in formula (37) p.13 of \cite{EMOT} shows that for all $u\neq 0, \nu\in (0, 1),$ one has
\begin{equation}
\label{fcos}
\int_0^{\infty}  \lbd^{\nu-1} \cos(\lbd u) \, d\lbd\; =\; \Ga(\nu) \cos (\pi\nu/2)\, \vert u\vert^{-\nu}. 
\end{equation}
The first statement of the lemma is hence simply a switching of the expectation and the integral. However, we cannot apply Fubini's theorem directly. Set $\mu$ for the probability distribution of $X.$ From (\ref{fcos}) and an integration by parts, we get
\begin{eqnarray*}
\Ga(\nu) \cos (\pi\nu/2)\,\E[ \vert X\vert^{-\nu}] & = & \int_\R \mu(du)\;\lpa\int_0^{\infty}  \lbd^{\nu-1} \cos(\lbd u) \,d\lbd\rpa\\
& = & (1-\nu) \int_\R \mu(du)\;\lpa \int_0^{\infty}  \frac{\sin(\lbd u)}{u}\, \lbd^{\nu-2}  \,d\lbd\rpa.
\end{eqnarray*}
Since
\begin{eqnarray*}
\int_\R \mu(du)\;\lpa \int_0^{\infty}  \lva \frac{\sin(\lbd u)}{u}\rva\, \lbd^{\nu-2}  \,d\lbd\rpa & \le & \int_\R \mu(du)\;\lpa \int_0^{\infty}  \lpa \lbd\wedge \frac{1}{\vert u\vert}\rpa\, \lbd^{\nu-2}  \,d\lambda\rpa\\
& \le & \frac{\E[ \vert X\vert^{-\nu}]}{\nu(1-\nu)}\; <\; +\infty,
\end{eqnarray*}
we may now apply Fubini's theorem and obtain
$$\Ga(\nu) \cos (\pi\nu/2)\,\E[ \vert X\vert^{-\nu}] \;= \; (1-\nu) \int_0^{\infty}  \lbd^{\nu-2}  \, d\lbd\; \lpa \int_\R \frac{\sin(\lbd u)}{u}\, \mu(du) \rpa.$$
The dominated convergence theorem entails that the function 
$$\psi(\lbd)\,=\, \int_\R \frac{\sin(\lbd u)}{u}\, \mu(du)$$
is differentiable, with derivative  
$$\psi^\prime(\lbd)\,=\,\int_\R \cos(\lbd u)\, \mu(du)\, =\, \E[\cos(\lbd X)].$$
Thus, another integration by parts yields 
$$\Ga(\nu) \cos (\pi\nu/2)\,\E[ \vert X\vert^{-\nu}]\; =\; \int_0^{\infty}  \lbd^{\nu-1}\, \E[\cos(\lbd X)] \, d\lbd\, -\,  \bigg[\lbd^{\nu-1}\psi(\lbd)\bigg]_0^{+\infty}$$
and it remains to prove that the bracket is zero. On the one hand, one has 
$$\lbd^{\nu -1}\vert\psi(\lbd)\vert \; \le \; \lbd^\nu\;\to\; 0\quad\mbox{as $\lbd \to 0$.}$$
On the other hand, using
$$\lbd^{\nu-1}\lva\frac{\sin (\lbd u)}{u}\rva\, \leq\, \lbd^{\nu -1}\frac{\vert \sin(\lbd u)\vert^{1-\nu}}{\vert u\vert}\, \leq \, |u|^{-\nu},$$
and the dominated convergence theorem, we see that $\lbd^{\nu -1}\vert\psi(\lbd)\vert \to 0$ as $\lbd \to +\infty.$ This completes the proof of the first statement of the lemma. The second statement may be handled similarly with the help of the formula
\begin{equation}
\label{fsin}
\int_0^{\infty}  \lbd^{\nu-1} \sin(\lbd x) \, d\lbd\; =\; \Ga(\nu) \sin (\pi\nu/2){\rm sgn}(x) \vert x\vert^{-\nu}, \qquad \vert \nu\vert < 1, 
\end{equation}
which is given e.g. in (38) p.13 in \cite{EMOT}.
 
\endproof

\begin{lemma} 
\label{FT} For all $x, y\in\R$ and $t \ge 0$ one has
$$\log(\E_{(x,y)}[e^{\i \lbd X_t}])\; =\; \i\lbd (x+ yt)\;-\;\frac{t^{\a +1}}{\a +1}\,(\i \lbd)^\a e^{-\i\pi\a\rho\, {\rm sgn}(\lbd)}, \qquad \lbd\in\R.$$
\end{lemma}

\proof It is clearly enough to consider the case $x = y =0.$ Integrating by parts yields the following representation of $X_t$ as a stable integral:
$$X_t\; =\; \int_0^\infty (t-s)^+ \,dL_s\; =\; \int_0^\infty (t-x)^+ \, M(dx),$$
where $M$ is an $\a-$stable random measure on $\R^+$ with Lebesgue control measure and constant skewness intensity $\b(x) =\b$ - see Example 3.3.3 in \cite{ST}. In the case $\a\neq 1,$ the statement of the lemma is a direct consequence of Proposition 3.4.1 (i) in \cite{ST}, reformulated with the $(\a,\rho)$ parametrization. In the case $\a =1,\rho =1/2$ we use Proposition 3.4.1 (ii) in \cite{ST} (with $\b =0$). The case $\a =1,\rho \neq 1/2$ follows from the symmetric case in adding a drift coefficient $\mu t$ for some $\mu\neq 0,$ which integrates in $\mu t^2/2.$

\endproof

We now set 
$$\sar\, = \,\frac{\sin (\pi\a(\rho -1/2))}{\a +1}\, \in \, (-1,1)\qquad\mbox{and}\qquad \car\, =\, \frac{\cos(\pi\a (\rho -1/2))}{\a +1}\, \in\, (0,1).$$
The proposition gives a representation for the Mellin transform of $X^+_t = X_t \Un_{\{X_t > 0\}}.$

\begin{proposition}
\label{lem:app1} For all $x, y\in\R, t > 0$ and $\nu \in (0,1)$ one has
$$\E_{(x,y)}[(X^+_t)^{-\nu}] = \frac{\Gamma(1-\nu)}{\pi}
 \int_0^{\infty} \! \lbd^{\nu-1} e^{-\car\lbd^\alpha t^{\alpha +1}} \sin(\lambda(x+yt) + \sar\lambda^\alpha t^{\alpha+1} + \pi\nu/2)\, d\lambda.$$
\end{proposition}
 
\proof Since $X_t$ is a stable random variable, it has a bounded density and $\E_{(x,y)}[(X^+_t)^{-\nu}]$ is hence finite for all $\nu\in (0,1).$ By Lemma \ref{FT} we have
$$\log(\E_{(x,y)}\left[e^{i\lambda X_t}\right])\, = \, \i\lambda (x + yt) - \lambda^\alpha t^{1+\alpha}(\car -\i s_{\a,\rho}), \quad \lambda\geq 0.$$
Taking the real part and integrating with respect to $\lambda^{\nu-1}$ on $]0,+\infty[$, we deduce
\begin{eqnarray*}
\int_0^{\infty}\! \lambda^{\nu-1} e^{- \car\lambda^\alpha t^{1+\alpha}} \cos(\lambda(x+yt) + \sar\lambda^\alpha t^{1+\alpha})\, d\lambda & = & \int_0^{\infty} \!\lambda^{\nu-1} \E_{(x,y)}\left[\cos(\lambda X_t)\right] \,d\lambda \\
& = & \Gamma(\nu) \cos\left(\frac{\pi\nu}{2}\right)\, \E_{(x,y)}\left[|X_t|^{-\nu}\right],
\end{eqnarray*}
where the second equality comes from Lemma \ref{switch}.
Similarly, taking the imaginary part entails
$$\int_0^{\infty}\! \lambda^{\nu-1} e^{- \car\lambda^\alpha t^{1+\alpha}} \sin(\lambda(x+yt) + \sar\lambda^\alpha t^{1+\alpha})\, d\lambda\;=\; \Gamma(\nu) \sin\left(\frac{\pi\nu}{2}\right)\, \E_{(x,y)}\left[|X_t|^{-\nu}\text{sgn}(X_t)\right].$$
Multiplying the first relation by $\sin(\pi\nu/2)$, the second by $\cos(\pi\nu/2),$ and summing, we finally obtain
$$\Gamma(\nu) \sin(\pi\nu)\,\E_{(x,y)}[(X^+_t)^{-\nu}] \,= 
\int_0^{\infty} \! \lbd^{\nu-1} e^{-\car\lbd^\alpha t^{\alpha +1}} \sin(\lambda(x+yt) + \sar\lambda^\alpha t^{\alpha+1} + \pi\nu/2)\, d\lambda,$$
which yields the required expression by the complement formula for the Gamma function.

\endproof
 
Our last proposition provides some crucial computations for the proof of Theorem B.

\begin{proposition} 
\label{lem:app2} Set $\nu\in (\a/(\a +1),1)$ and $s = (1-\nu)(\a+1)\in (0,1).$

\medskip

{\em (i)} For every $y > 0,$ one has
$$\int_0^\infty \E_{(0,y)} [(X^+_t)^{-\nu}]\, dt\; =\; (\a+1)^{1-\nu}\Ga(1-s) \sin(\pi s (1-\g))\frac{\Ga(1-\nu)^2 }{\pi}\,y^{s-1}\cdot$$

{\em (ii)} For every $y < 0,$ one has
$$\int_0^\infty \E_{(0,y)} [(X^+_t)^{-\nu}]\, dt\; =\;(\a+1)^{1-\nu}\Ga(1-s)\sin(\pi\g s) \frac{\Ga(1-\nu)^2 }{\pi}\, \vert y\vert^{s-1}\cdot$$

{\em (iii)} For every $x<0,$ one has
$$\int_0^\infty \E_{(x,0)} [(X^+_t)^{-\nu}]\, dt\; =\;(\a +1)^{-\frac{\alpha}{\alpha+1}} \,\Gamma\left(\frac{1-s}{\a +1}\right)\sin (\pi\g)\Ga\left(\frac{1}{\a+1}\right) \frac{\Ga(1-\nu) }{\pi}\,  \vert x\vert^{\frac{s-1}{\a +1}}.$$
\end{proposition}

\proof Suppose first $x = 0$ and $y\in\R.$ Integrating the expression on the right-hand side of Proposition \ref{lem:app1} yields a double integral of the form
\begin{eqnarray*} 
I_\nu & = &\int_0^{\infty}\lpa\int_0^{\infty} \lambda^{\nu-1} e^{- \car\lambda^\alpha t^{1+\alpha}} \sin(\lambda yt + \sar\lambda^\alpha t^{1+\alpha} +\nu\pi/2)\, d\lambda\rpa dt\\
\quad\qquad&= &\int_0^{\infty}\lpa\int_0^{\infty} r^{\nu-1} e^{-\car r^\a} \sin(ryt^{-1/\alpha}  + \sar r^\alpha+\nu\pi/2)\, dr\rpa t^{-\nu (1+1/\alpha)}\, dt\\
\quad\qquad&= &\int_0^{\infty}\lpa\int_0^{\infty} r^{\nu-1} e^{-\car r^\a}  \sin(ryu + \sar r^\alpha+\nu\pi/2)\, dr\rpa \a u^{-s}\, du\\
\quad\qquad&= &\int_0^{\infty}\lpa\int_0^{\infty} \a u^{-s} \sin(ryu + \sar r^\alpha+\nu\pi/2)\, du\rpa r^{\nu-1} e^{-\car r^\a} \, dr,
\end{eqnarray*}
where the first, resp. second, equality comes from the change of variable $\lambda t^{1+1/\alpha}=r$, resp. $u= t^{-1/\alpha}$, and the switching of the integrals in the third equality is made exactly as in Lemma \ref{switch}, using the fact that $s\in (0,1)$ and $s+\nu > 1.$ \\

Suppose first $y > 0.$ We start by computing the integral in $u$ with the help of formul\ae\,\! (\ref{fcos}) and (\ref{fsin}) and some trigonometry:
$$\a\int_0^{\infty} u^{-s} \sin(ryu + \sar r^\alpha+\nu\pi/2)\, du
\;=\; \a\Gamma(1-s) \cos ((s-\nu)\pi/2 - \sar r^\a )  (yr)^{s-1}.$$
We then compute the integral in $r$ with the change of variable $z=r^\alpha,$ using the notation $Z = e^{\i\pi\a(\rho - 1/2)}$:
\begin{eqnarray*}
I_\nu & = & \a\Gamma(1-s) \,y^{s-1}\int_0^{\infty} r^{\a(1-\nu)-1} e^{-\car r^\a} \cos ((s-\nu)\pi/2 - \sar r^\a) \, dr\\
&= &\Gamma(1-s)\,y^{s-1} \int_0^{\infty} z^{-\nu} e^{-\car z} \cos((s-\nu)\pi/2  - \sar z ) \, dz\\
& = & (1+\a)^{1-\nu}\Ga(1-s)\Ga(1-\nu)\Re( e^{\i\pi(s-\nu)/2} Z^{\nu -1})\,y^{s-1}\\
& = & (1+\a)^{1-\nu}\Ga(1-s)\Ga(1-\nu)\sin(\pi s (1-\g))\,y^{s-1},
\end{eqnarray*}
where the third line follows after some algebraic simplifications. By Proposition \ref{lem:app1}, this completes the proof of (i).\\

Suppose now $y < 0.$ An analogous computation to the above shows that
$$I_\nu\; =\;\a\int_0^{\infty} u^{-s} \sin(ryu + \sar r^\alpha+\nu\pi/2)\, du
\;=\; \a\Gamma(1-s) \sin ((s+\nu-1)\pi/2 + \sar r^\a )  \vert yr\vert^{s-1}.$$
The integral in $r$ is then computed in the same way and yields the formula
\begin{eqnarray*}
I_\nu& =&(1+\a)^{1-\nu}\Ga(1-s)\Ga(1-\nu)\,\Im( e^{\i\pi(s+\nu-1)/2} {\bar Z}^{\nu -1})\,\vert y\vert ^{s-1}\\
 & = & (1+\a)^{1-\nu}\Ga(1-s)\Ga(1-\nu)\sin(\pi \g s)\,\vert y\vert ^{s-1},
 \end{eqnarray*}
which completes the proof of (ii) by Proposition \ref{lem:app1}.\\

We last suppose $x < 0$ and $y =0.$ We again integrate the expression on the right-hand side of Proposition \ref{lem:app1}, making the changes of variable $\lambda t^{1+1/\alpha}=r$ and $u= t^{-(1+1/\alpha)}.$ This yields a double integral of the form\\
$$ \frac{\a}{\a+1}\int_0^{\infty}\lpa\int_0^{\infty} r^{\nu-1} e^{-\car r^\a}  \sin(rxu + \sar r^\alpha+\nu\pi/2)\, dr\rpa u^{-(s+\a)/(1+\a)}\, du,$$
where we can switch the orders of integration as in Lemma \ref{switch} because $(s+\a)/(1+\a)\in (0,1)$ and $(s+\a)/(1+\a) + \nu > 1.$ We then compute the integral in $u$ similarly as above and find
$$\frac{\a}{\a+1} \,\Gamma\left(\frac{1-s}{\a +1}\right) \,\sin (\pi\a/2(\a +1) + \sar r^\a ) \, \vert xr\vert^{\frac{s-1}{\a +1}}.$$
We finally compute the integral in $r$ with the change of variable $r = z^{1/\a},$ and get after some algebraic manipulations
$$I_\nu\; =\;(\a +1)^{-\frac{\alpha}{\alpha+1}}\sin (\pi\g)\Ga\left(\frac{1}{\a+1}\right) \,\Gamma\left(\frac{1-s}{\a +1}\right)  \vert x\vert^{\frac{s-1}{\a +1}},$$
which completes the proof of (iii) by Proposition \ref{lem:app1}. 

\endproof

\begin{remark} It seems hard to find an explicit formula in general for
$$\int_0^\infty \E_{(x,y)} [(X^+_t)^{-\nu}]\, dt$$
when $(x,y)$ is not on a coordinate axis. In the symmetric Cauchy case, some further computations show that the integral equals
$$\frac{1}{\sin (\pi\nu)}\,\Im\lpa \int_0^\infty (-(x+yt +\i t^2/2))^{-\nu}\,dt\rpa.$$
This can be rewritten with the hypergeometric function, apparently not in a tractable manner when $x y \neq 0.$
\end{remark}

\section{Proof of Theorem B}

The following lemma shows the aforementioned and intuitively obvious fact that $T_0$ is a proper random variable for any starting point.

\begin{lemma}
\label{Klar}
For all $x, y\in \R$ one has $\Pb_{(x,y)}[ T_0 < +\infty] = 1.$
\end{lemma}

\proof Suppose first $x =-1$ and $y=0.$ Then
$$\Pb_{(-1,0)}[ T_0 = +\infty]\, =\, \Pb_{(0,0)}[ X_\infty^* < 1] \,= \,\Pb_{(0,0)}[ X_\infty^* =0] \,\le\, \Pb_{(0,0)}[ X_1 \le 0] \,<\, 1,$$  
where the second equality comes from the self-similarity of $X$ and the strict inequality from the fact that $X_1$ is a two-sided stable random variable - see Lemma \ref{FT}. On the other hand, setting $T = \inf\{ t > 0, X_t > 0\},$ it is clear by self-similarity that under the probability measure $\Pb_{(0,0)}$ one has
$$T\, \elaw\, k T$$
for all $k > 0.$  In particular, $\Pb_{(0,0)}[T\in \{0, +\infty\}] = 1.$ Moreover, the zero-one law for the Markov process $(X,L)$ entails that $\Pb_{(0,0)}[T =0]$ is $0$ or $1.$ Since $\Pb_{(0,0)}[T=+\infty]  = \Pb_{(0,0)}[ X_\infty^* =0] < 1,$ we get $\Pb_{(0,0)}[T=+\infty] = 0$ whence $\Pb_{(-1,0)}[ T_0 = +\infty] = 0$ as desired. Notice that it also entails $\Pb_{(0,0)}[T=0] = 1,$ as mentioned in the introduction.

Using again self-similarity, this entails $\Pb_{(x,0)}[ T_0 <+\infty] = 1$ for all $x\le 0,$ and also for all $x \ge 0$ in considering the dual process $-L.$ The fact that $\Pb_{(x,y)}[ T_0 <+\infty] = 1$ for all $x,y$ such that $xy < 0$ follows then by a comparison of the sample paths.

Suppose now that $x \le 0, y <0.$ Introduce the stopping time $S = \inf\{t > 0, L_t > 0\},$ which is finite a.s. under $\Pb_{(x,y)}$ because $\vert L\vert$ is not a subordinator. It is clear that $L_S \ge 0$ and $X_S < 0$ a.s. Applying the strong Markov property, we see from the above cases that
$$\Pb_{(x,y)}[ T_0 = +\infty] \, \le \, \Pb_{(x,y)}[ \Pb_{(X_S,L_S)}[ T_0 = +\infty]] = 0.$$
The same argument holds for $x\ge 0, y > 0.$

\endproof

Assume now $x < 0$ or $x=0$ and $y < 0.$ It is clear that at $T_0$ the process $X$ has a non-negative speed, which entails by right-continuity that $L_{T_0} \ge  0$ a.s. Applying the Markov property at $T_0$ entails
\begin{equation}
\label{Mark}
\pbxy [X_t \in du]\; =\; \int_0^\infty \int_0^t \Pb_{(0,z)}[ X_{t-s} \in du]\,\pbxy [T_0\in ds, L_{T_0}\in dz]
\end{equation}
for all $t, u > 0.$ Integrating in time yields then after a change of variable and Fubini's theorem
$$\int_0^\infty \pbxy [X_t \in du]\, dt\; =\; \int_0^\infty \lpa\int_0^\infty \Pb_{(0,z)}[ X_t\in du]\,dt\rpa\pbxy [L_{T_0}\in dz]$$
for all $u > 0.$ Integrating in space along $u^{-\nu}$ and applying again Fubini's theorem shows finally the general formula
\begin{eqnarray}
\label{Gen}
\int_0^\infty \E_{(x,y)} [(X^+_t)^{-\nu}] dt & =& \int_0^\infty \Pb_{(x,y)}[L_{T_0} \in dz] \lpa\int_0^\infty \E_{(0,z)} [(X^+_t)^{-\nu}] dt\rpa
\end{eqnarray}
which is valid for all $\nu\in\R,$ with possibly infinite values on both sides.

\subsection{Proof of (i)} Assume $x =0$ and $y < 0.$ Setting $\nu \in (\a/(\a+1), 1),$ a straightforward application of Proposition \ref{lem:app2} (i) and (ii) is that both sides of (\ref{Gen}) are finite, which leads to
$$\E_{(0,y)} [L_{T_0}^{s-1}]\; =\; \vert y\vert^{s-1}\lpa \frac{\sin (\pi \g s)}{\sin (\pi (1-\g)s)}\rpa$$
for all $s\in (0,1).$ The formula extends then to $\{\vert s \vert < 1/(1-\g)\}$ by analytic continuation. On the other hand, for all $\mu\in (0,1)$ and $s\in (-1,1),$ the formula
$$\int_0^\infty \frac{\sin(\pi\mu) x^s}{\pi\mu (x^2 + 2 \cos (\pi\mu) x + 1)}\, dx\; =\; \frac{\sin (\pi\mu s)}{\mu \sin (\pi s)}$$
is a simple and well-known consequence of the residue theorem. Recalling that
$$\chi\, =\, \frac{\g}{1-\g}\, \in\, (0,1)$$
and the definition of $\CC_\mu,$ we deduce
$$\E_{(0,y)} [L_{T_0}^{s-1}]\; =\; \vert y\vert^{s-1}\E[\CC_\chi^{(1-\g)s}]$$
for all $\vert s \vert < 1/(1-\g),$ which concludes the proof of (i) by Mellin inversion.

\qed

\subsection{Proof of (ii)} Assume $x < 0$ and $y =0.$ Another application of (\ref{Gen}) combined with Proposition \ref{lem:app2} (i) and (iii) shows that
\begin{equation}
\label{Melvil}
\E_{(x,0)}[L_{T_0}^{s-1}]\; =\;  \frac{(1+\alpha)^{\frac{1-s}{1+\alpha}}\Ga(\frac{\a +2}{\a +1})\Ga (\frac{1-s}{\a +1})\sin (\pi\g)}{\Ga(\frac{s}{\a +1})\Ga (1-s)\sin (\pi s(1-\g))}\,\vert x\vert^{\frac{s-1}{\a+1}}
\end{equation}
for all $s\in (0,1).$ A simple analysis on the Gamma factors shows that the above expression remains finite for all $\vert s \vert < 1/(1-\g).$

\qed

\subsection{Some further Mellin inversions} In this paragraph we would like to invert (\ref{Melvil}) for certain values of the parametrization $(\a, \rho).$ Without loss of generality we set $x = -1, y =0.$ Applying the complement formula for the Gamma function, we first deduce from (\ref{Melvil})
\begin{equation}
\label{Poupaud}
\E_{(-1,0)}[L_{T_0}^{s-1}]\; =\;  (1+\alpha)^{\frac{1-s}{1+\alpha}}\;\frac{\Ga(\frac{\a +2}{\a +1})\Ga (\frac{1-s}{\a +1})\Ga (s(1-\g))\Ga(1-s(1-\g))}{\Ga(\frac{s}{\a +1})\Ga (1-s)\Ga(\g)\Ga(1-\g)}
\end{equation}
for $\vert s \vert < 1/(1-\g).$

\subsubsection{The Cauchy case} We have $\a = 1$ and $\rho \in (0,1),$ whence $\g =\rho/2\in (0,1/2).$ As mentioned in the introduction, set  
$$\delta \; =\; \frac{1}{2(1-\g)}\; =\; \frac{1+\chi}{2}\; \in (1/2, 1).$$
Applying the Legendre-Gauss multiplication formula transforms (\ref{Poupaud}) into
$$\E_{(-1,0)}[L_{T_0}^{s-1}]\; =\; 2^{\frac{s-1}{2}}\;\times\; \frac{\sin(\pi\g)\sin(\pi s/2)}{\sin(\pi s(1-\g))}\cdot$$
As above, this entails that under $\Pb_{(x,0)}$ one has
$$L_{T_0}\, \elaw\, \sqrt{2\vert x\vert}\, (\CC_\delta^{1-\g})^{(1)},$$
which provides an striking similarity with the law of $L_{T_0}$ under $\Pb_{(0,y)}$ for $y < 0.$ Notice that these two laws are however never the same, because $\delta \neq \chi.$ 

\subsubsection{The Brownian case} We have $\a =2, \rho = \chi = 1/2$ and $\g = 1/3.$ Applying three times the Legendre-Gauss multiplication formula and simplifying the quotients shows
$$\E_{(-1,0)}[L_{T_0}^{s-1}]\; =\;  9^{\frac{s-1}{3}}\;\times\;\frac{\Ga(1/2 +s/3)}{\Ga (5/6)}\;\times\;\frac{\Ga(1/2 -s/3)\Ga(1/3)}{\Ga(2/3 -s/3)\Ga (1/6)}$$
for all $s\in (0,1).$ Inverting the Mellin transform, this entails that under $\Pb_{(x,0)}$ one has
$$L_{T_0}\, \elaw\, \vert 9x\vert^{1/3}\lpa \frac{\GG_{5/6}}{\BB_{1/6,1/6}}\rpa^{1/3},$$
where $\GG_c$ resp. $\BB_{a,b}$ stands for the standard Gamma resp. Beta random variable, and the quotient is assumed independent. Gor'kov \cite{Gk} provides an expression of the density of $L_{T_0}$ under $\pbxy$ in terms of the confluent hypergeometric function - see also formula (3) in \cite{La}. It seems however that the above simple identity in law has passed unnoticed in the literature on integrated Brownian motion. 

\begin{remark} It is well-known that $\log(\GG_c)$ and $\log(\BB_{a,b})$ are infinitely divisible random variables, and this property is hence also shared by $\log(L_{T_0})$ under $\Pb_{(x,0)}.$ The question whether $L_{T_0}$ itself is infinitely divisible is an interesting open problem for Brownian motion.
\end{remark}

\subsubsection{The case $\a < 1$} We have $\rho \in (0,1),\g \in (0,1/2)$ and $\chi\in (0,1).$ Set  
$$\eta \; =\; \frac{1}{(\a +1)(1-\g)}\; =\; \frac{1}{1+\alpha(1-\rho)}\; \in (1/2, 1)\qquad\mbox{and}\qquad \sigma\; =\; \frac{\a +1}{2}\in (1/2,1).$$
To express our result, we need some further notation. For every $\mu\in(0,1)$ set $\ZZ_\mu$ for the standard positive $\mu-$stable random variable \cite{Z}, which is characterized through its Mellin transformation by
$$\E[\ZZ_\mu^{s}] \, =\, \frac{\Ga(1-\frac{s}{\mu})}{\Ga(1-s)}, \qquad s < \mu.$$
Applying again the Legendre-Gauss formula entails
\begin{eqnarray*}
\E_{(-1,0)}[L_{T_0}^{s-1}] & =&  \kappa \,\lpa\frac{2}{(1+\a)^{\frac{1}{1+\a}}}\rpa^{s-1}\times\;\frac{\Ga(1 +s(1-\g))}{\Ga (1+\frac{s}{\a+1})}\;\times\;\frac{\Ga(1-s(1-\g))}{\Ga(1 -\frac{s}{2})}\;\times\;\frac{\Ga(1+\frac{1-s}{\a+1})}{\Ga (1+ \frac{1-s}{2})}\\
& = & \kappa \,\lpa\frac{2}{(1+\a)^{\frac{1}{1+\a}}}\rpa^{s-1}\E[\ZZ_\eta^{-\frac{s}{\a+1}}]\,\times\, \E[\ZZ_\delta^{\frac{s}{2}}]\,\times\, \E[\ZZ_\sigma^{\frac{s-1}{2}}]
\end{eqnarray*}
for all $s\in (0,1),$ where $\kappa$ is the normalizing constant. Identifying, this shows that under $\Pb_{(x,0)}$ one has
$$L_{T_0}\, \elaw\, 2\,\lva \frac{x}{\a+1}\rva^{\frac{1}{\a+1}}\,\ZZ_\sigma^{\frac{1}{2}}\,\times\lpa \frac{\ZZ_\delta^{\frac{1}{2}}}{\ZZ_\eta^{\frac{1}{\a+1}}}\rpa^{(1)},$$
where the product and the quotient are assumed independent.

\begin{remark} The above argument shows that the function
$$\M_\a(s)\; =\;\frac{\Ga (\frac{1}{\a+1}+\frac{s-1}{\a+1})\Ga(1 +\frac{s-1}{\a+1})}{\Ga (\frac{1}{\a+1})\Ga(1 +(s-1))}$$
is the Mellin transform of a positive random variable for all $\a\le 1.$ This is equivalent to the fact that the independent product $\X_u  = \GG_u^u\times \GG_1^u$ is an exponential mixture for all $u =1/(\a+1)\in [1/2,1).$ It is easy to see that $\X_u$ is also an exponential mixture for all $u \ge 1.$ However, this property is not true in general for $u < 1/2.$ Taking $u = 1/3$ viz. $\a = 2,$ we see indeed from the Legendre-Gauss multiplication formula that
$$\M_2(s)\; =\;\lpa\frac{4}{27}\rpa^{\frac{s-1}{3}}\,\frac{\Ga(\frac{2}{3})}{\Ga(\frac{2}{3} +\frac{s-1}{3})}$$
is log-concave and hence not the Mellin transform of a positive measure.

\end{remark}

\subsubsection{The case $1<\a < 2$} We first separate (\ref{Poupaud}) into

\begin{eqnarray*}
\E_{(-1,0)}[L_{T_0}^{s-1}] & = &  (1+\alpha)^{\frac{1-s}{1+\alpha}}\times\;\frac{\Ga(\frac{1}{\a +1})\Ga (s(1-\g))}{\Ga(1-\g)\Ga(\frac{s}{\a +1})}\,\times\,\frac{(\frac{1}{\a +1})\Ga (\frac{1-s}{\a +1})\Ga(1-s(1-\g))}{\Ga(\g)\Ga (1-s)}\\
& =& (1+\alpha)^{\frac{1-s}{1+\alpha}}\times\;\E\lcr \ZZ_{1/(1+\a(1-\rho))}^{-s/(\a+1)}\rcr\,\times\,\frac{(\frac{1}{\a +1})\Ga (\frac{1-s}{\a +1})\Ga(1-s(1-\g))}{\Ga(\g)\Ga (1-s)}\cdot
\end{eqnarray*}
We next make the following assumption
$$\g\; \le \;1/3.$$
Notice that this assumption is fulfilled in the spectrally negative case, where $\rho = 1-1/\a$ viz. $\g = (\a-1)/(\a +1) < 1/3.$ Setting $\M_{\a,\g}(s)$ for the second multiplicand on the right-hand side, we use again the Legendre-Gauss multiplication formula to get the transformation
\begin{eqnarray*}
\M_{\a,\g}(s)& =& \kappa \,3^s \;\frac{\Ga (\frac{1-s}{\a +1})}{\Ga (\frac{1-s}{3})}\,\times\,\frac{\Ga (\frac{1}{2}-\frac{s}{3})}{\Ga (\frac{2}{3}-\frac{s}{3})}\,\times\,\frac{\Ga(1-s(1-\g))}{\Ga (\frac{1}{2}-\frac{s}{3})\Ga (1-\frac{s}{3})}\\
& =& {\tilde \kappa} \,\lpa 3.2^{-2/3}\rpa^s\, \E\lcr \ZZ_{\frac{\a+1}{3}}^{\frac{s-1}{3}}\rcr\,\times\,\E\lcr \BB_{1/6,1/6}^{\frac{1-s}{3}}\rcr\,\times\,\frac{\Ga(1-s(1-\g))}{\Ga (1-\frac{2s}{3})}\\
& =& {\tilde \kappa} \,\lpa 3.2^{-2/3}\rpa^s\, \E\lcr \ZZ_{\frac{\a+1}{3}}^{\frac{s-1}{3}}\rcr\,\times\,\E\lcr \BB_{1/6,1/6}^{\frac{1-s}{3}}\rcr\,\times\,\E\lcr \ZZ_{2/3(1-\g)}^{2s/3}\rcr
\end{eqnarray*}
where $\kappa, {\tilde \kappa}$ are normalizing constants. This leads to the identity in law 
$$L_{T_0}\; \elaw\; 3. 2^{-2/3}\lva \frac{x}{\a+1}\rva^{\frac{1}{\a+1}}\,\times\,\lpa \frac{\ZZ_{\frac{\a+1}{3}}}{\BB_{1/6,1/6}}\rpa^{\frac{1}{3}}\times\;\lpa \frac{\ZZ_{2/3(1-\g)}^{2/3}}{\ZZ_{1/(1+\a(1-\rho))}^{1/(\a+1)}}\rpa^{(1)},$$
an extension of the Brownian case because when $\a =2$ the first multiplicand is $\BB_{1/6,1/6}^{-1/3},$ whereas the second one reads
$$\lpa \ZZ_{1/2}^{-1/3}\rpa^{(1)} \elaw\; 2^{2/3} \lpa \GG_{1/2}^{1/3}\rpa^{(1)}\elaw\; 2^{2/3} \GG_{5/6}^{1/3}.$$

\begin{remark} We do not know whether $\M_{\a,\g}(s)$ is still the Mellin transform of a positive random variable in the remaining case $1/3<\g\le 1/(\a+1).$ This would be a consequence of the exponential mixture property of the independent product $\GG_u^u\times\GG_1^{1-u}$ for all $u\in (1/3, 1/2).$ It is easily shown that the latter property holds for all $u=1/n, n\ge 2,$ and we think that it does for all $u\in [0, 1/2].$ 
\end{remark}

\section{Proof of Theorem A}

We first reduce the problem to the situation where the bivariate process $(X,L)$ starts from a coordinate axis.

\begin{lemma}
\label{Gene}
Assume $x < 0.$ For all $y\in\R$ one has
$$\pbxy[T_0 > t]\; \asymp\;  \Pb_{(x, 0)}[T_0 > t], \qquad t\to +\infty.$$
\end{lemma}

\proof Fix $t > 1$ and suppose first $y > 0.$ One has $\pbxy[T_0 > t]\le\Pb_{(x,0)}[T_0 > t]$ by a direct comparison of the sample paths. On the other hand,
\begin{eqnarray*}
\pbxy[T_0 > t] & \ge & \pbxy [X_1 < x, L_1 < 0, T_0 > t]\\
& = & \pbxy \lcr X_1 < x, X_1^* < 0, L_1 < 0, \,\Pb_{(X_1, L_1)}[T_0 > t-1]\rcr\\
& \ge & \pbxy [X_1 < x, X_1^* < 0, L_1 < 0]\,\times\,\Pb_{(x, 0)}[T_0 > t-1]\;\; \ge \; c\, \Pb_{(x, 0)}[T_0 > t]
\end{eqnarray*}
for some $c >0,$ where the equality follows from the Markov property, the second inequality from a comparison of the sample paths, and the third inequality from a support theorem in uniform norm for the L\'evy stable process $L$. More precisely, the latter process can be approximated by any continuous function in uniform norm, because the support of its L\'evy measure is $\R$ - see Corollary 1 in \cite{TS0}. 

\medskip

Fix again $t > 1$ and suppose now $y < 0.$ Then $\pbxy[T_0 > t]\ge\Pb_{(x,0)}[T_0 > t],$ and similarly as above one has
\begin{eqnarray*}
 \Pb_{(x, 0)}[T_0 > t] & \ge &  \Pb_{(x, 0)} \lcr X_1 < x, X_1^* < 0, L_1 < y, \,\Pb_{(X_1, L_1)}[T_0 > t-1]\rcr\\
& \ge & \Pb_{(x, 0)} [X_1 < x, X_1^* < 0, L_1 < y]\,\times\,\pbxy[T_0 > t-1]\;\;\ge \; c\, \pbxy[T_0 > t]
\end{eqnarray*}
for some $c > 0.$ This completes the proof.

\endproof

In the remainder of this section,  we will  implicitly assume, without loss of generality, that  
$$\{x=0, y <0\}\qquad\mbox{or}\qquad\{x < 0, y =0\}.$$ 
We start by studying the asymptotics at infinity of the density function of $L_{T_0}$ under $\pbxy,$ which we denote by $\foxy.$

\begin{lemma}\label{lem:asympLT0} There exists $c > 0$ such that
$$\foxy(z)\; \sim\; c z^{-1/(1-\g)}, \qquad z \to +\infty.$$
\end{lemma}

\proof
If $x =0,$ the asymptotic is a direct consequence of the explicit expression of $f^0_{(0,y)}$ which is given in Theorem B (i). If $y =0,$ Theorem B (ii) shows that the first positive pole of the Mellin transform of $L_{T_0}$ under $\Pb_{(x,0)}$ is at $1/(1-\g)$, and is simple. The required asymptotic for $f^0_{(x,0)}$ is then a consequence of a converse mapping theorem for Mellin transforms - see e.g. Theorem 6.4 in \cite{Ja}.

\endproof

\begin{remark} (a) The converse mapping theorem for Mellin transforms yields also an explicit expression for the constant $c$, but we shall not need this information in the sequel.

\medskip

(b) We believe that the above asymptotic remains true for $x < 0$ and all $y\neq 0.$ However, the Mellin transform of $L_{T_0}$ under $\pbxy$ is then expressed with the help of a double integral which is absolutely divergent, and whose singularities are  difficult to study at first sight.

\medskip

(c) The lemma entails by integration that 
$$\pbxy[L_{T_0} > z]\; \sim \; c\chi^{-1}\, z^{-\chi}, \qquad z \to +\infty.$$
Heuristically, it is tempting to write by scaling $L_{T_0} = T_0^{1/\a} \vert L_1\vert$ and since $\pbxy[ \vert L_1\vert > z] \sim c z^{-\a} \ll z^{-\chi}$ at infinity, we may infer that
$$\pbxy [T_0 > t] \; \asymp\; t^{-\chi/\a} \; =\; t^{-\theta}, \qquad t \to +\infty.$$
This explains the equivalence between finite moments stated in the introduction. We will prove in the remainder of this section that this heuristic is actually correct.
\end{remark}

The following lemma provides our key-estimate.

\begin{lemma}
\label{lem:H} For all $\nu\in (\a(1-\theta)/(\a+1), 1)$ there exists $c > 0$ such that 
$$\E_{(x,y)}\left[ \int_0^t \Un_{\{T_0>t-u\}}\,\E_{(0,L_{T_0})}\left[ (X_u^+)^{-\nu}\right]  du \right]\; \sim\;c\, t^{1-(1+1/\alpha)\nu - \theta},\qquad t\to +\infty.$$
\end{lemma}

\proof We first assume $\nu\in (\a/(\a+1), 1)$ and transform the expression on the left-hand side. From (\ref{Mark}), Fubini's theorem, and the Markov property, we obtain
$$\int_0^{\infty} e^{-\lambda t}\,\E_{(x,y)}\left[(X_t^+)^{-\nu}\right]dt
\;=\;\E_{(x,y)}\left[e^{-\lambda T_0}  \int_0^{\infty} e^{-\lambda t}\,\E_{(0,L_{T_0})}\left[(X_t^+)^{-\nu}\right]dt \right], \quad \lbd \ge 0,$$
both sides being finite because of Proposition \ref{lem:app2}. Integrating by parts shows then, with the help of (\ref{Gen}) and Proposition \ref{lem:app2}, that
\begin{multline*}
\lambda \int_0^{\infty} e^{-\lambda t}\int_t^{\infty} \left(\E_{(x,y)}\left[ (X_u^+)^{-\nu}\right]- \E_{(x,y)}\left[ \E_{(0,L_{T_0})}\left[ (X_u^+)^{-\nu}\right] \right] \right)du\, dt
\end{multline*}
\vspace*{-.3cm}
\begin{eqnarray*}
&= &\E_{(x,y)}\left[  (1-e^{-\lambda T_0})\int_0^{\infty} e^{-\lambda t}\E_{(0,L_{T_0})}\left[ (X_t^+)^{-\nu}\right] dt \right]\\
& = &\E_{(x,y)}\left[ \int_0^{\infty} \lbd \,e^{-\lambda t}\lpa \int_0^t \Un_{\{T_0>t-u\}} \E_{(0,L_{T_0})}\left[ (X_t^+)^{-\nu}\right] du\rpa dt \right].
\end{eqnarray*}
Inverting the Laplace transforms shows that
\begin{equation}\label{eq:invLap}
\E_{(x,y)}\left[ \int_0^t \Un_{\{T_0>t-u\}}\E_{(0,L_{T_0})}\left[ (X_u^+)^{-\nu}\right]  du \right] \; = \; H_{(x,y)}(t),
\end{equation}
with the notation
\begin{equation}
\label{Hxy}
H_{(x,y)}(t)\;=\;\int_t^{+\infty} \left(\E_{(x,y)}\left[ (X_u^+)^{-\nu}\right]- \E_{(x,y)}\left[ \E_{(0,L_{T_0})}\left[ (X_u^+)^{-\nu}\right] \right] \right)du, \qquad t >0.
\end{equation}
It remains therefore to compute the asymptotics of the function $H_{(x,y)}$, which only depends on the law of $L_{T_0}$ under $\pbxy$. From Proposition \ref{lem:app1}, the additive property of sine, and a change of variable, we get

\medskip

${\displaystyle
\E_{(x,y)}\left[ (X_u^+)^{-\nu}\right]- \E_{(x,y)}\left[\E_{(0,L_{T_0})}\left[ (X_u^+)^{-\nu}\right]\right]}$
\begin{eqnarray}
\notag&= &\frac{2\Gamma(1-\nu)}{\pi}\int_0^{\infty} \lambda^{\nu-1} e^{-c_{\alpha, \rho} \lambda^\alpha u^{\alpha+1}} \Phi_u(\lambda u^{1+1/\alpha}) \,d\lambda\\
\label{eq:E-E}& = &\frac{2\Gamma(1-\nu)}{\pi}\frac{1}{u^{(1+1/\alpha)\nu}}\int_0^{\infty} \xi^{\nu-1} e^{-c_{\alpha, \rho} \xi^\alpha }  \Phi_u(\xi)\, d\xi 
\end{eqnarray}
where the function $\Phi_u$ is defined by 
\begin{eqnarray*}
\Phi_u(\xi) &= &\E_{(x,y)}\left[\cos\left(\frac{\xi x}{2 u^{1+1/\alpha}}+s_{\alpha,\rho}\,\xi^\alpha +\frac{\nu\pi}{2}+ \frac{\xi(y+L_{T_0})}{2u^{1/\alpha}}\right)  \sin\left(\frac{\xi x}{2 u^{1+1/\alpha}}+ \frac{\xi(y-L_{T_0})}{2u^{1/\alpha}}\right)\right]\\
&= &\!\!\!\int_0^{\infty} \!\!\!\cos\left(\frac{\xi x}{2 u^{1+1/\alpha}}+s_{\alpha,\rho}\,\xi^\alpha +\frac{\nu\pi}{2}+ \frac{\xi(y+z)}{2u^{1/\alpha}}\right)  \sin\left(\frac{\xi x}{2 u^{1+1/\alpha}}+ \frac{\xi(y-z)}{2u^{1/\alpha}}\right)\foxy (z) \,dz. 
\end{eqnarray*}
Setting $F_u(\xi,z)$ for the trigonometric function inside the integral, a change of variable entails
\begin{equation}
\label{chg}
u^\theta\lpa\int_0^\infty \! F_u(\xi,z)\,\foxy (z) \,dz\rpa\; =\; \int_0^\infty \! F_u(\xi,ru^{1/\a})\,u^{\frac{1}{\a(1-\g)}}\foxy (ru^{1/\a}) \,dr.
\end{equation}
A further application of the converse mapping theorem and of Theorem B shows that
$\foxy (z) \sim c z^{1/(1-\g)}$ as $z \to 0+,$ for some $c > 0.$ This estimate and Lemma \ref{lem:asympLT0} yield the uniform bound
$$\foxy(z)\; \le \;K z^{-1/(1-\g)}, \qquad z > 0$$
for some $K > 0.$ Hence, for $u$ large enough, the integrated function on the right-hand side of (\ref{chg}) is dominated by $K(r\wedge 1)\, r^{-(\chi +1)}$ 
for some $K > 0,$ which is an integrable function because $\chi \in (0,1).$ Applying the dominated convergence theorem in (\ref{chg}) and using Lemma \ref{lem:asympLT0} shows that there exists $c > 0$ such that for all $\xi > 0$ 
\begin{eqnarray*}
u^\theta\,\Phi_u(\xi)& \rightarrow & - c \int_0^{\infty}  \cos( s_{\alpha, \rho}\, \xi^\alpha +(\nu\pi +\xi r)/2)  \sin(\xi r/2)r^{-(1+\chi)}\, dr\\
&= & -(c/2) \Gamma(-\chi) \sin( s_{\alpha, \rho} \xi^\alpha +(\nu -\chi)\pi/2))\, \xi^{\chi}
\end{eqnarray*}
as $u\to \infty,$ where the equality comes from (\ref{fsin}), an integration by parts and some trigonometry. 
Plugging back this expression in (\ref{eq:E-E}), we deduce by dominated convergence that
\begin{multline}
\label{glou}
\E_{(x,y)}\left[ (X_u^+)^{-\nu}\right]- \E_{(x,y)}\left[\E_{(0,L_{T_0})}\left[ (X_u^+)^{-\nu}\right]\right] \\
\!\!\!\!\!\!\!\!\equi_{u \rightarrow +\infty} c \lpa\int_0^{\infty} \xi^{\chi +\nu-1} e^{-c_{\alpha, \rho} \xi^\alpha }\sin( s_{\alpha, \rho} \xi^\alpha +(\nu -\chi)\pi/2)) d\xi\rpa u^{-((1+1/\a)\nu +\theta)}
\end{multline}
for some $c > 0.$ A last computation shows that the integral on the right-hand side equals
$$(\a+1)^{\frac{\chi +\nu}{\a}}\,\Ga\left(\frac{\chi +\nu}{\a}\right)\,\sin (\pi(\rho\nu+ (\rho -1)\chi))$$
and is positive because $\rho\nu+ (\rho -1)\chi\in (0,1)$
for all $\nu\in (\a(1-\theta)/(\a +1), 1)$ as can be readily checked. The final result follows then by integration, and the proof is complete for $\nu\in (\a/(\a +1), 1).$

\medskip

Suppose now $\nu\in (\a(1-\theta)/(\a +1), \a/(\a +1)).$ The left-hand side of (\ref{eq:invLap}) is well-defined and the estimate (\ref{glou}), which does not require the lower bound $\nu > \a/(\a +1),$ together with the positivity of the constant entails that the integral in (\ref{Hxy}) is absolutely convergent, because  
$(1+1/\a)\nu +\theta > 1.$ By analytic continuation this shows that (\ref{eq:invLap}) remains valid for $\nu\in (\a(1-\theta)/(\a +1), \a/(\a +1)),$ and the estimate (\ref{glou})
holds as well. This completes the proof, again by integration of (\ref{glou}).

\endproof

\subsection{Proof of the upper bound} Fix $A>0$ and $\nu\in (\a/(\a +1), 1).$ By continuity and positivity there exists $\eps>0$ such that for all $z\in[0,A]$,
$$\int_0^1  \E_{(0,z)}\left[ (X_u^+)^{-\nu}   \right] \,du \;\geq\; \eps.$$
For all $t > 0,$ we get from (\ref{eq:invLap}), a change of variable and the self-similarity
\begin{align*}
\notag t^{(1+1/\alpha)\nu+\theta -1} H_{(x,y)}(t) & \geq \; t^{(1+1/\alpha)\nu+\theta -1}\,\E_{(x,y)}\left[ \Un_{\{T_0>t\}}\int_0^t  \E_{(0,L_{T_0})}\left[ (X_u^+)^{-\nu}\right]  du\right] \\
\notag & =\;t^{\theta}\,  \E_{(x,y)}\left[ \Un_{\{T_0>t\}}\int_0^1  \E_{(0,\frac{1}{t^{1/\alpha}}L_{T_0})}\left[ (X_u^+)^{-\nu}\right]  du\right]\\ 
\notag &\geq \;\eps t^{\theta}\, \pbxy[ T_0>t, L_{T_0}\leq A t^{1/\alpha}]  \\
\notag &\geq \;\eps t^{\theta}\, \left(\Pb_{(x,y)}[T_0>t]  - \Pb_{(x,y)}[T_0>t, L_{T_0}\geq A t^{1/\alpha}] \right)\\
\label{eq:H>}&\geq \;\eps t^{\theta}\, \left(\Pb_{(x,y)}[T_0>t]  - \Pb_{(x,y)}[T_0>t]^{1-1/p} \Pb_{(x,y)}[L_{T_0}\geq A t^{1/\alpha}]^{1/p} \right),
\end{align*}
where the last inequality follows from H\"older's inequality and is valid for all $p>1$. We now take $t> 1$ and 
$$\frac{1}{p}\; = \;1-\frac{1}{\ln(t)}\cdot$$
On the one hand, Lemma \ref{lem:asympLT0} entails 
\begin{multline*}
\limsup_{t\rightarrow +\infty}\;  t^{\theta}\,\Pb_{(x,y)}[T_0>t]^{\frac{1}{\ln(t)}} \,\Pb_{(x,y)}[L_{T_0}\geq A t^{1/\alpha}]^{1-\frac{1}{\ln(t)}}\\
 \leq \; \limsup_{t\rightarrow +\infty} \; t^{\theta}\,\Pb_{(x,y)}[L_{T_0}\geq A t^{1/\alpha}]^{1-\frac{1}{\ln(t)}}\;=\; K \;<\; +\infty.
\end{multline*}
On the other hand, Lemma \ref{lem:H} shows that
$$ t^{(1+1/\alpha)\nu+\theta -1} H_{(x,y)}(t)\; \to \; c \;>\;  0\qquad \mbox{as $t\to +\infty.$}$$
Putting everything together entails
$$t^\theta\Pb_{(x,y)}[T_0>t]\; \le \; {\tilde K}$$
for some finite ${\tilde K}$ as soon as $t$ is large enough.

\qed

\subsection{Proof of the lower bound}

We start with the following lemma :

\begin{lemma}\label{lem:intPb}
One has 
$$\int_0^t \Pb_{(x,y)}[T_0>u]\, du\; \asymp \; t^{1-\theta}\qquad \mbox{as $t\rightarrow +\infty.$}$$
\end{lemma}

\proof
Firstly, integrating the above upper bound for $\Pb_{(x,y)}[T_0>t]$ entails the existence of a finite $\kappa_2$ such that
$$\int_0^t \Pb_{(x,y)}[T_0>u]\, du\;\le\; \kappa_2 \,t^{1-\theta}\qquad \mbox{as $t\rightarrow +\infty.$}$$
To prove the lower inequality, we fix $\nu \in (\alpha(1-\theta)/(1+\alpha),\alpha/(1+\alpha))$ and deduce from Proposition \ref{lem:app1} the uniform bound
\begin{equation}
\label{UB}
\E_{(0,y)}[(X^+_u)^{-\nu}] \;\leq\;  \frac{\Gamma(1-\nu)}{\pi}
 \int_0^{\infty} \! \lbd^{\nu-1} e^{-\car\lbd^\alpha u^{\alpha +1}} \, d\lambda \;\leq \;Ku^{-\nu(1+1/\alpha)}, \qquad u >0,
\end{equation}
for some finite constant $K.$ Set $\eta =\nu(1+1/\alpha) \in (0,1)$ and fix $\eps\in (0,1).$ Using (\ref{eq:invLap}) and (\ref{UB}) we decompose
\begin{eqnarray*}
t^{\eta +\theta -1} H_{(x,y)}(t)& \leq & K t^{\eta +\theta -1}\lpa \int_0^{t(1-\varepsilon)}  \frac{\Pb_{(x,y)}[T_0>u] }{(t-u)^{\eta}}\,du \;+\; \int_{t(1-\varepsilon)}^t   \frac{\Pb_{(x,y)}[T_0>u] }{(t-u)^{\eta}}\,du\rpa \\
&\leq & K \eps^{-\eta} t^{\theta -1}\int_0^{t} \Pb_{(x,y)}[T_0>u] \,du \;+ \;\frac{K t^\theta\varepsilon^{1-\eta}}{1-\eta} \,\Pb_{(x,y)}[T_0>t(1-\varepsilon)]\\  
&\leq & {\tilde K} \eps^{-\eta} \lpa t^{\theta -1}\int_0^{t} \Pb_{(x,y)}[T_0>u] \,du \;+ \;\varepsilon\rpa
\end{eqnarray*}
for some finite ${\tilde K}$, where the third inequality follows from the upper bound. Applying Lemma \ref{lem:H} and taking $\varepsilon$ small enough shows finally that there exists $\kappa_1> 0$ such that 
$$\int_0^t \Pb_{(x,y)}[T_0>u]\, du\;\ge\; \kappa_1\, t^{1-\theta}\qquad \mbox{as $t\rightarrow +\infty.$}$$
\endproof

\noindent
We can now finish the proof. Fixing $A > 0$ and applying the mean value theorem entails 
$$A\, t^\theta\, \Pb_{(x,y)}[T_0>t] \; \ge \; t^{\theta -1} \int_t^{t+tA}   \Pb_{(x,y)}[T_0>u] \,du\; \ge \; \kappa_1 (1+A)^{1-\theta} \, -\, \kappa_2$$
as $t\rightarrow +\infty$, for some constants $0 < \kappa_1 < \kappa_2 < \infty$ given by Lemma \ref{lem:intPb}. Since $\theta <1,$ the lower bound follows in choosing $A$ large enough.

\qed

\medskip
  
\noindent
{\bf Acknowledgement.}  Ce travail a b\'en\'efici\'e d'une aide de la Chaire {\em March\'es en Mutation}, F\'ed\'eration Bancaire Fran\c{c}aise.

\end{document}